\documentclass[12pt]{amsart}
\usepackage{amsfonts,amssymb,amsthm,amsmath,amsxtra,amscd,verbatim,eucal}
\usepackage[all]{xy}
\usepackage[dvips]{graphics}

\theoremstyle{plain}
\newtheorem{thm}{Theorem}[section]
\newtheorem{lem}[thm]{Lemma}

\theoremstyle{definition}

\theoremstyle{remark}

\newtheorem{rmk}[thm]{Remark}

\newtheorem{question}[thm]{Question}

\def\N{{\mathbb N}}

\def\A{{\mathbb A}}
\def\R{{\mathbb R}}
\def\Q{{\mathbb Q}}
\def\P{{\mathbb P}}

\def\O{\mathcal{O}}

\def\M{\mathcal{M}}

\def\a{\mathfrak{a}}
\def\bb{\mathfrak{b}}
\def\cc{\mathfrak{c}}

\def\ee{{\bf e}}
\def\cc{{\bf c}}                           

\def\aa{\alpha}

\def\m{\mu}

\def\t{\tau}

\def\o{\circ}
\def\ov{\overline}

\def\.{\cdot}
\def\({\Big{(}}
\def\){\Big{)}}
\def\~{\widetilde}
\def\^{\widehat}
\def\-{}

\DeclareMathOperator{\codim} {codim}

\DeclareMathOperator{\Spec} {Spec}

\DeclareMathOperator{\Bl} {Bl}

\DeclareMathOperator{\mult} {mult}

\DeclareMathOperator{\lc}{lc}

\DeclareMathOperator{\Supp} {Supp}

\DeclareMathOperator{\Vol} {Vol}

\begin{document}

\title{Multiplicities and log canonical threshold}

\author[T. de Fernex]{Tommaso de Fernex}
\address{Department of Mathematics, Statistics and
Computer Science, University of Illinois at Chicago,
851 Morgan St., M/C. 249, Chicago, IL 60607-7045,USA}
\email{{\tt defernex@math.uic.edu}}

\author[L. Ein]{Lawrence Ein}
\address{Department of Mathematics, Statistics and
Computer Science, University of Illinois at Chicago,
851 Morgan St., M/C. 249, Chicago, IL 60607-7045, USA}
\email{{\tt ein@math.uic.edu}}

\author[M. Musta\c{t}\v{a}]{Mircea~Musta\c{t}\v{a}}
\address{Isaac Newton Institute for Mathematical Sciences,
20 Clarkson Road, Cambridge, CB3 0EH, England}
\email{{\tt mirceamustata@yahoo.com}}

\thanks{Research of the second author was partially supported by 
NSF~Grant DMS~99-70295. The third author
served as a Clay Mathematics Institute Long-Term Prize Fellow
while this research has been done.}

\maketitle

\section*{Introduction}

Let $X$ be a smooth algebraic variety, defined over an algebraically closed
field, of characteristic zero. Consider a closed subscheme $Y\subset X$, and
let $Z$ be an irreducible component of $Y$, of codimension $n$.
 We study the connection between
two of the basic invariants of the ideal $I$ defining $Y$ in $\O_{X,Z}$.

The first invariant is very classical: it is the Samuel multiplicity $e(I)$
of $\O_{X,Z}$ along $I$. The second invariant has been much studied recently,
due to its proeminent role in higher dimensional birational geometry.
It is called the log canonical threshold of $\O_{X,Z}$ along $I$,
and it is denoted by $\lc(I)$.

Our main result is the following:

\begin{thm}\label{thm1-intro}
With the above notation, we have
\begin{equation}\label{eq_main}
e(I)\geq\frac{n^n}{\lc(I)^n}.
\end{equation}
\end{thm}

The case $n=2$ of the above theorem 
was proved by Corti in \cite{corti}. We also characterize the case when in
(\ref{eq_main}) we have equality: this happens if and only if
the integral closure of $I$ is a power of the maximal ideal in $\O_{X,Z}$.

If $Z\subseteq Y$, but it is not an irreducible component of $Y$,
then $I$ is not zero-dimensional, so $e(I)$ is not defined. In this case,
Gaffney and Gassler have introduced in \cite{GG} the Segre numbers
of $\O_{X,Z}$ along $I$, denoted by $e_i(I)$, for $i=1,\ldots,n$.
One can formulate a possible generalization of Theorem~\ref{thm1-intro}
to this setting (see Question~\ref{open_question} below). 
We prove this in the case
$n=2$. More precisely, we show the following:

\begin{thm}\label{thm2-intro}
If $Z\subset X$ is a closed and irreducible subset of codimension $2$,
contained in the support of $Y$, and if $e_1(I)$ and $e_2(I)$
are the Segre numbers of $\O_{X,Z}$ along $I$, then
\begin{equation}\label{eq_surface}
\frac{4}{\lc(I)}e_1(I)+e_2(I)\geq\frac{4}{\lc(I)^2}.
\end{equation}
\end{thm}

In fact, in the two-dimensional case, we give also a stronger version
of this inequality (see Theorem~\ref{thm5} below).

The proofs of these results are based on deformation to monomial ideals.
To fix the ideas, we discuss now the case of Theorem~\ref{thm1-intro}.
We prove first an inequality between the length and the log canonical 
threshold:
\begin{equation}\label{eq_main2}
l(R/I)>\frac{n^n}{n!\cdot\lc(I)^n}.
\end{equation}
The inequality in Theorem~\ref{thm1-intro} follows by applying
(\ref{eq_main2}) to the powers $I^t$, dividing by $t^n$, and passing to 
the limit.

In order to prove~(\ref{eq_main2}),
standard arguments allow us to reduce the general case to that
of an ideal $J$
in a polynomial ring. By considering a deformation to a monomial ideal
${\rm in}(J)$, it is enough to prove~(\ref{eq_main2}) when $J$ is monomial.
We achieve this via a combinatorial argument.

Interest in bounds for log canonical thresholds is partly motivated
by applications to birational rigidity. Iskovskikh and Manin proved
in \cite{IM} that a smooth quartic in ${\mathbb P}^4$
is birationally rigid (hence, in particular, it
is not rational), and this has started a whole direction in higher
dimensional birational geometry (see, for example \cite{Pu2},
\cite{Pu3}, \cite{CPR}). Recently, Corti initiated in \cite{corti}
a new approach to this circle of ideas, based on Shokurov's Connectedness
principle and 
the two-dimensional case of inequality (\ref{eq_main}).
In this spirit, we 
apply Theorem~\ref{thm1-intro}
in \cite{DEM} to give 
an upper bound 
for the log canonical threshold of a generic
projection of a projective variety. This improves a bound due to Pukhlikov
from \cite{Pu1}, bound which plays a major role in his proof of the fact
that a smooth hypersurface of degree $n$ in ${\mathbb P}^n$,
 $n\geq 4$, is birationally
rigid.

\subsection{Acknowledgements}
We are grateful to Rob Lazarsfeld for many useful conversations.

\section{Lower-bounds to length and multiplicity}

Throughout this paper, all varieties are assumed to be defined over
an algebraically closed field $k$ of characteristic zero. We first
fix the notation.

Let $X$ be a smooth variety, $Y\subset X$ a closed subscheme, and
let $Z$ be an irreducible component of $Y$. We denote by
$I\subset\O_{X,Z}$ the ideal locally defining $Y$. Let $\lc(I)$
be the log canonical threshold of $(\O_{X,Z},I)$, and set
$\mu(I)=1/\lc(I)$. For the definition and basic properties of the
log canonical threshold, we refer to \cite{kollar}. Note that
following the definition in \cite{kollar},
$\lc(I)$ is the maximum of $\lc(U,U\cap Y)$, when $U\subseteq X$ 
is an open subset such that $U\cap Z\neq\emptyset$.

Let $\codim_XZ=n$. We denote by $l(\O_{X,Z}/I)$ the length of
$\O_{X,Z}/I$, and by $e(I)$ the multiplicity of $\O_{X,Z}$ along $I$.
Our main result is the following theorem.

\begin{thm}\label{thm1}
With the above notation, we have
\begin{equation}\label{thm1-eq}
l(\O_{X,Z}/I) > \frac{n^n \m(I)^n}{n!}.
\end{equation}
\end{thm}

The above result easily gives the theorem we have stated in the Introduction.

\begin{thm}\label{cor1}
Keeping the above notation, we have
\begin{equation}\label{cor1-eq}
e(I) \geq n^n \m(I)^n.
\end{equation}
\end{thm}

\begin{proof}
By applying Theorem~\ref{thm1} to powers $I^t$ of $I$,
we get 
\begin{equation}\label{l(I^t)}
l(\O_{X,Z}/I^t) \geq \frac{n^n \m(I^t)^n}{n!}.
\end{equation}
where $\m(I^t)$ is the inverse of the log canonical threshold of 
$(\O_{X,Z},I^t)$. Since
$$
\m(I^t) = t \m(I) \quad\text{and}\quad 
\lim_t \frac{n!\cdot l(\O_{X,Z}/I^t)}{t^n} = e(I),
$$
(\ref{cor1-eq}) follows by multiplying both sides of~(\ref{l(I^t)})
by $n!/t^n$ and passing to the limit as $t \to \infty$.
\end{proof}

The proof of Theorem~\ref{thm1} consists of a reduction to 
the case of monomial ideals, when the assertion follows from the
explicit description of the log canonical threshold.
Before proving Theorem~\ref{thm1}, 
we recall the definition and the main properties of the Newton polytope 
associated to a monomial ideal.
To a monomial $x_1^{a_1} \cdots x_n^{a_n}$ we associate the point
$(a_1,\dots,a_n) \in \N^n \subset \R_+^n$ (here 
$\R_+ = [0,+\infty)$). If $J \subset R = K[x_1,\dots,x_n]$ 
is a monomial ideal, the Newton polytope
$P(J)$ associated to $J$ is defined as the convex hull in $\R_+^n$ 
of the points corresponding to the monomials in $J$. 
Let $\lc(J)$ be the log canonical threshold of $(R,J)$, 
and set $\mu(J) = 1/\lc(J)$.
Then a result due to Howald~\cite{howald} implies that
$$
\mu(J) = \min \{ \aa > 0 \;|\; \aa \cdot \ee \in P(J) \},
$$ 
where $\ee = (1,\dots,1) \in \R_+^n$. 
We will also write $\mu(P(J))$ for $\mu(J)$. Moreover, we
extend this notation, defining
$$
\mu(P) = \inf \{ \aa > 0 \;|\; \aa \cdot \ee \in P \},
$$ 
for any subset $P \subset \R_+^n$, such that $\alpha\cdot\ee\in P$, for some
$\alpha>0$.

We observe that $J$ is a zero
dimensional ideal if and only if $\R_+^n\setminus P(J)$ is
bounded in $\R_+^n$. We define the 
volume of a zero dimensional monomial ideal $J$ by setting
$$
\Vol(J) := n! \Vol(\R_+^n\setminus P(J)),
$$
where the right hand side is the volume of the corresponding bounded region,
in the Euclidean metric. 
\footnote{For a zero dimensional monomial ideal $J$, $\Vol(J) = e(J)$,
but we will not need this fact.}

\begin{lem}\label{lem2}
If $J \subset R$ is a zero dimensional monomial ideal, then
$n!\cdot l(R/J) > \Vol(J)$.
\end{lem}

\begin{proof} 
We denote by $(u_1,\dots,u_n)$ the coordinates on $\R^n$.
Let 
$$
C_0 = \{ (u_1,\dots,u_n) \in \R_+^n \;|\; 0 \leq u_i \leq 1 \}
$$
be the standard unitary $n$-cube of $\R_+^n$.
For any point $a = (a_1,\dots,a_n) \in (\N^n\setminus P(J))$, let
$\t_a : \R_+^n \to \R_+^n$ be the translation
which sends the origin to the point $a$, 
and consider the $n$-cube $C_a = \t_a(C_0)$.
Then 
\begin{align*}
l(R/J) &= \# \{\,\text{monomials}\, \not \in J \} \\
&\geq \#(\N^n\setminus P(J)) = 
\Vol( \bigcup_{a \in (\N^n\setminus P(J))} C_a).
\end{align*}
Note that a point $u = (u_1,\dots,u_n) \in (\R_+^n\setminus P(J))$
is contained in the $n$-cube $C_{[u]}$,
where $[u] = ([u_1],\dots,[u_n]) \in (\N^n\setminus P(J))$
($[\;\,]$ denotes the integral part operator).
Therefore the union of the $n$-cubes appearing on the right hand side
covers $\R^n\setminus P(J)$. This proves the weak inequality in the statement.
To prove that the inequality is strict, it is enough to note that 
$\R_+^n\setminus \bigcup_{a \in (\N^n\setminus P(J))} C_a$ is not convex,
hence 
$\bigcup_{a \in (\N^n\setminus P(J))} C_a \cap P(J)$ has positive volume.
\end{proof}

\begin{proof}[Proof of Theorem~\ref{thm1}]
By passing to completion, 
we obtain a zero dimensional ideal $\^ I \subset \^ \O_{X,Z}$.
Note that $\^ \O_{X,Z}$ is isomorphic to a formal power
series ring $K[[x_1,\ldots,x_n]]$.
We will identify these two 
$k$-algebras
 under some fixed isomorphism. Then, since $\^ I$ is zero dimensional,
we can find an ideal $J$ in the polynomial ring $R = K[x_1,\dots,x_n]$, 
with $\Supp(V(J)) = \{0\} \subset
 \A^n$, and whose completion in $K[[x_1,\dots,x_n]]$
is $\^ I$. 

Since we have an isomorphism of $k$-algebras $\O_{X,Z}\simeq R/J$,
it is clear that $l(\O_{X,Z}/I)=l(R/J)$. Moreover, we have
$\lc(\O_{X,Z},I)=\lc(R,J)$. To see this, recall that 
if $W\subset \widetilde{W}$ is a subscheme of a smooth variety
$\widetilde{W}$
 over some field $F$, then $\dim(\widetilde{W})-\lc(\widetilde{W},W)$
depends only on $W$, but not on the embedding $W\subset\widetilde{W}$
(see, for example, \cite{mustata}, Corollary~3.5). Therefore it is enough
to note that computing $\lc(\A_K^n, V(J))$ (over $K$) is the same as computing
$\lc(\A_k^n\times T, V(J))$ (over $k$), where $T$ is a suitable smooth
variety over $k$, with function field $K$. 

The above argument shows that we may reduce ourselves to the case of 
the ideal $J\subset R$. By an extension of scalars to the algebraic
closure of $K$, we may also assume that $K$ is algebraically closed.
The proof now splits into two parts. First we prove the theorem
assuming that $J$ is a monomial ideal. Then we deduce the inequality
for arbitrary $J$ by a degeneration argument.

Suppose that $J$ is a monomial ideal. Let $F$ be the facet of $P(J)$
containing the point $\m \. \ee$. Since $F$ is bounded (recall that 
$J$ is zero dimensional), the equation of the hyperplane supporting
$F$ can be written as
$$
\sum_{i=1}^n u_i/a_i = 1,
$$
for some $a_i > 0$. This gives
$$
\Vol(J) \geq \prod_{i=1}^n a_i \quad\text{and}\quad
\m(J) = \left(\sum_{i=1}^n \frac{1}{a_i}\right)^{-1}.
$$
Then, by comparing the arithmetic and geometric means of $\{1/a_i\}_i$, we 
see that
$$
\Vol(J) \geq n^n \mu(J)^n.
$$
Therefore the theorem follows, in the monomial case, by
Lemma~\ref{lem2}.

Now we consider an arbitrary ideal $J \subset R = K[x_1,\dots,x_n]$
with $\Supp(R/J) =\{0\} \subset \A^n$.
After fixing a multiplicative order on the coordinates $x_i$,
we can take a deformation to the initial monomial ideal (see
\cite{eisenbud}, Chapter 15).
This is a flat family  $\{J_s\}_{s\in K}$ such that $S/J_s \cong S/J$,
for all $s \ne 0$, and such that ${\rm in}(J):=J_0$ is a monomial ideal. 
Since this is a flat deformation, $\Supp(R/{\rm in}(J)) = \{0\}$ and 
$l(R/J) = l(R/{\rm in}(J))$. Moreover, the semicontinuity property
of the log canonical threshold 
(see, for example, \cite{DK}) gives $\lc(R,J) \geq \lc(R,{\rm in}(J))$,
hence $\m(J) \leq \m({\rm in}(J))$. Thus inequality~(\ref{cor1-eq}) follows 
now from the case of monomial ideals.
\end{proof}

The boundary case in Theorem~\ref{cor1} is characterized in the
following theorem. 

\begin{thm}\label{cor1=}
Under the assumptions of Theorem~\ref{cor1},
\begin{equation}\label{cor1=-eq}
e(I) = n^n \m(I)^n
\end{equation}
if and only if there is a positive integer $q$, such that
the integral closure $\ov I$ of $I$ is equal to
 $\M^q$, where $\M$ is the maximal ideal of $\O_{X,Z}$. 
Moreover, in this case $q=n\m(I)$.
\end{thm}

\begin{proof}
If we assume that $\ov I = \M^q$ for some positive integer $q$, 
then  
$e(I)=q^n$ and $\lc(I)=n/q$, hence~(\ref{cor1=-eq}) 
is satisfied.

Conversely, assume that~(\ref{cor1=-eq}) holds. 
Since $e(I) \in \N$ and $\m(I) \in \Q_+$, we
see that $n\m(I) \in \N$. We will prove that
$$
\ov I = \M^{n \m(I)}.
$$
By hypothesis, $e(I) = e(\M^{n \m(I)})$, and we will show that 
\begin{equation}\label{inclusion}
I \subseteq \M^{n \m(I)}.
\end{equation}
This will imply
our assertion by a theorem of Rees \cite{Rees}, which says
that if two zero dimensional ideals $I\subseteq J\subset \O_{X,Z}$
are such that $e(I)=e(J)$, then $\overline{I}=\overline{J}$.

Therefore it is enough to prove~(\ref{inclusion}). Arguing as in the 
beginning of the proof of Theorem~\ref{thm1}, we reduce
to study the case when $I$ is an ideal in $R = K[x_1,\dots,x_n]$
with $\Supp(R/I) =\{0\} \subset \A^n$.

We deform all the powers $I^t$ to monomial ideals ${\rm in}(I^t)$,
as follows. Consider a fixed multiplicative order on the monomials in $R$.
We first deform $I^t$ to the tangent cone, i.e.,
to the ideal generated by the sum of terms of lowest degree in $f$,
for $f\in I^t$, and then take the initial ideal of the resulting ideal,
with respect to the monomial order. 
By the way we have made the deformation, $I\subseteq\M^p$
if and only if ${\rm in}(I)\subseteq\M^p$. Note also that
${\rm in}(I^t)\cdot {\rm in}(I^s)\subseteq {\rm in}(I^{s+t})$.

As in the proof of Theorem~\ref{thm1}, for every 
positive integer $t$,
$$
\frac{n!\cdot l(R/{\rm in}(I^t))}{t^n} \geq \frac{n^n \mu({\rm in}
(I^t))^n}{t^n}
\geq \frac{n^n \mu(I^t)^n}{t^n} = n^n \mu(I)^n.
$$
On the other hand, since we are assuming that $e(I) = n^n \mu(I)^n$, 
we also have
$$
\frac{n!\cdot l(R/{\rm in}(I^t))}{t^n} = \frac{n!\cdot l(R/I^t)}{t^n} \to
e(I) = n^n \mu(I)^n,
$$
for $t \to \infty$. Combining these two formulas, we obtain
$$
\lim_{t\to\infty} \frac{\mu({\rm in}(I^t))}t = \mu(I).
$$
Set 
$$
P_t = \frac 1t P({\rm in}(I^t)) \quad\text{and}\quad 
P_{\infty} = \cup_{r \in \N} P_{2^r}.
$$
Note that for every pair of positive integers $r \leq s$, 
$P_{2^r} \subseteq P_{2^s}$. This follows from the fact that
$$
({\rm in}(I^{2^r}))^{2^{s-r}} \subseteq {\rm in}(I^{2^s}).
$$
Therefore $P_{\infty}$ is a convex subset in $\R_+^n$, and 
$$
\m(P_{\infty}) = \lim_{r\to\infty} \frac{\m({\rm in}(I^{2^r}))}{2^r} = \m(I).
$$
On the other hand, an application of Lemma~\ref{lem2} gives
\begin{align*}
n! \Vol(\R_+^n\setminus P_{\infty}) 
&= \lim_{r\to\infty} (n! \Vol(\R_+^n\setminus P_{2^r})) \\
&\leq \lim_{r\to\infty} \frac{n!\cdot l(R/{\rm in}(I^{2^r}))}{(2^r)^n} 
= n^n \m(I)^n.
\end{align*}
By the convexity of $P_{\infty}$, we can find an hyperplane 
$F \subset \R_+^n$ passing through the point $\m(P_{\infty}) \cdot \ee$
and disjoint from the interior of $P_{\infty}$.
Let $\sum_{i=1}^n u_i/a_i = 1$ be the equation of $F$, and denote
by $S_F$ the simplex in $\R_+^n$ containing the origin and 
having $F$ as diagonal facet. Then we have
$$
n^n \m(I)^n \leq n! \Vol(S_F) \leq n! \Vol(\R_+^n\setminus 
P_{\infty}) = n^n \m(I)^n.
$$
Thus this is a chain of equalities. 
This implies the equality between the arithmetic and geometric means
of $\{1/a_i\}_i$, which can happen only if all $a_i$ are equal.
Therefore we conclude that $a_i = n\m(I)$ for all $i$. 
Observing that $P({\rm in}(I)) = P_1 \subseteq P_{\infty}$, we deduce that
${\rm in}(I)\subseteq\M^{n\m(I)}$.
Therefore $I \subseteq \M^{n\m(I)}$, which concludes the proof of the theorem.
\end{proof}

In the remaining part of this section, we discuss a more general set up.
As before, we consider a closed subscheme $Y$ of a smooth algebraic
variety $X$ together with 
an irreducible subvariety $Z$ of $X$. However, now we only 
assume that $Z$ is contained in $Y$. In other words, 
the ideal $I \subset \O_{X,Z}$, locally defining $Y$, is not 
necessarily zero dimensional.
We would like to generalize inequality~(\ref{cor1-eq})
in this setting. 

Segre numbers, introduced by 
Gaffney and Gassler in~\cite{GG}, appear as
a natural choice for a substitute for the Samuel multiplicity
in this more general context. Like the 
Samuel multiplicity of a zero dimensional ideal, the Segre numbers 
of the ideal $I$ can be computed as intersection numbers
and have a natural interpretation as multiplicities, related
to the Vogel cycle associated to $I$.
The following seems a plausible generalization of Theorem~\ref{cor1}.

\begin{question}\label{open_question}
With the above notation, if $n=\codim_XZ$, and if
$e_1(I),\ldots,e_n(I)$ are the Segre numbers of $I$, is it true that 
\begin{equation}\label{eq_open}
\sum_{k=1}^n \frac{e_k(I)}{k^k \m(I)^k}  \geq 1 ?
\end{equation}
\end{question}

Note that if $I\subset\O_{X,Z}$ is a complete intersection ideal,
of dimension $r$, then $e_i(I)=0$ for $i\neq n-r$, and
$e_{n-r}(I)=e(I+(f_1,\ldots,f_r))$, where $f_1,\ldots,f_r$
are general linear combinations of a system of generators of $I$.
In particular, the inequality~(\ref{eq_open}) holds,
by applying Theorem~\ref{cor1=} for the ring $\O_{X,Z}/(f_1,\ldots,f_r)$.
This also shows that the coefficients before the Segre numbers in
(\ref{eq_open}) are optimal. 

In the next section, we will give a positive answer to the
above question in the codimension two case.

\section{Inequalities in codimension two}

In this section we concentrate on the codimension two case. 
Specifically, we show that the inequality in
Question~\ref{open_question} holds in this case. We 
provide also a sharper inequality when the ambient variety 
is a surface. First we prove the following lemma.

\begin{lem}\label{lem1}
Let $I$ be a monomial ideal in the ring $R = K[x_1,x_2]$, and let
$\m(I)$ denote the inverse of the log canonical threshold $\lc(I)$ of $(R,I)$.
Write $I = x_1^{b_1}x_2^{b_2} \. \a$, where $\a$ is a zero dimensional ideal.
Then the length of $R/\a$ is bounded by:
\begin{equation}\label{lem1-eq}
l(R/\a) \geq 2 (\m(I) - b_1)(\m(I) - b_2).
\end{equation}
\end{lem}

\begin{proof}
Note first that 
$$
\m(I) \geq \max \{ b_1, b_2 \}.
$$
If we have equality, then (\ref{lem1-eq}) is 
trivially satisfied, so from now on we assume
that this is not the case.

Let $P(I)$ and $P(\a)$ be 
the Newton polytopes associated to $I$ and $\a$, respectively, and let
$$
U = \{ (u_1,u_2) \in \R_+^2 \;|\; u_i \geq b_i \}.
$$
If $\t_b$ is the translation sending the origin
of $\R_+^2$ to the point $b = (b_1,b_2)$, we have
$$
\t_b(\R_+^2\setminus P(\a)) = U\setminus P(I).
$$
The boundary of $P(I)$ is supported on the union of 
finitely many lines, among which the two lines 
$l_1 : u_1 = b_1$ and $l_2 : u_2 = b_2$.
Let $l$ be the line supporting a facet of $P(I)$ which contains
the point $(\m(I),\m(I))$. Note that $l$ is neither $l_1$ nor $l_2$,
since we have assumed $\mu(I)>\max\{b_1,b_2\}$. 
Let $T \subset \R_+^2$ be the triangular region bounded by 
the three lines $l, l_1, l_2$.
Because of the convexity of $P(I)$, we see that $T \subseteq U\setminus P(I)$.
Then we have
$$
\Vol(\a) \geq 2\Vol(T),
$$
hence, by Lemma~\ref{lem2},
$$
l(R/\a) \geq \Vol(T).
$$
Let $S \subset \R_+^2$ be the square with vertices $(0,0)$, $(\m(I),0)$,
$(0,\m(I))$, and $(\m(I),\m(I))$.
$S \cap U$ is a rectangle of area $(\m(I) - b_1)(\m(I) - b_2)$, and 
$$
S \cap U \subset T.
$$
Then inequality~(\ref{lem1-eq}) follows 
by observing that the area of a rectangle 
inscribed in a right triangle (with the edges of the rectangle 
parallel to the legs of the triangle)
does not exceed half of the area of the triangle.
This is an easy consequence of the inequality
between the arithmetic and the geometric means. (Alternatively, one can give 
a synthetic proof of this fact by simply drawing a picture
and suitably doubling the rectangle.)
\end{proof}

Before applying Lemma~\ref{lem1}, we fix some notation.
Let $X$ be a smooth variety, 
$Z \subset X$ an irreducible subvariety of codimension 2, and
$Y \subset X$ a subscheme containing $Z$.
Let $I \subset \O_{X,Z}$ be the ideal locally defining $Y$.
We can write $I = f \. \a$, where $f \in \O_{X,Z}$ and $\a$ is zero 
dimensional. We denote by $\m(I)$ the inverse of the log canonical 
threshold of $(R,I)$, by $\mult_Z(f)$ the multiplicity of $f$
at the generic point of $Z$ and by $e(\a)$ the Samuel multiplicity
of $\a$. We have the following theorem.

\begin{thm}\label{thm4}
With the above notation, we have 
\begin{equation}\label{thm4-eq}
4 \m(I) \mult_Z(f) + e(\a) \geq 4 \m(I)^2. 
\end{equation}
\end{thm}

\begin{proof}
Note that it is enough to prove that the following inequality holds:
\begin{equation}\label{eq4}
4\m(I)\mult_Z(f) + 2l(\O_{X,Z}/\a) \geq 4\m(I)^2.
\end{equation}
Indeed, since $I^t = f^t \cdot \a^t$ for every $t \in \N$,
we can apply (\ref{eq4}) to powers $I^t$ of $I$, so that
inequality~(\ref{thm4-eq}) follows by dividing both sides by 
$t^2$ and passing to the limit as $t \to \infty$, as in the proof
of Theorem~\ref{cor1}.

The idea is to deform $I$ to a monomial ideal, and to deduce the
inequality~(\ref{eq4}) from Lemma~\ref{lem1}. If the log canonical center
of $(\O_{X,Z},\lc(I)\cdot I)$ 
has codimension one, then we are done,
as in this case $\mu(I)\leq\mult_Z(f)$. Therefore we may assume that this
is not the case.

We proceed as in the proof
of Theorem~\ref{thm1}.
Passing to completion, we obtain an ideal $\^ I$ in 
$\^ \O_{X,Z}$. We identify $\^\O_{X,Z}$ with $K[[x_1,x_2]]$
via a fixed isomorphism. This gives a decomposition
$\^ I=g\cdot\bb$, with $\bb$ a zero dimensional ideal.
Since $\bb$ is zero dimensional, there is an ideal
$\cc\subset R=K[x_1,x_2]$ such that $\bb=\^\cc$. 

We fix now a monomial order on $R$ and then consider 
${\rm in}(g)$ and ${\rm in}(\cc)$. Here ${\rm in}(g)=x_1^{b_1}x_2^{b_2}$
is the largest of the monomials in $g$ of smallest degree, and 
${\rm in}(\cc)$ is the initial ideal, under the fixed order, of the ideal
obtained from $\cc$ after deforming to the normal cone. We put $J=x_1^{b_1}
x_2^{b_2}\cdot {\rm in}(\cc)$.
As before, we have 
$$\lc(I)=\lc(g\cdot\cc)\geq\lc(J),$$
and $l(\O_{X,Z}/I)=l(R/x_1^{b_1}x_2^{b_2}\cdot {\rm in}(\cc))$.
Note also that $\mult_Z(f)=b_1+b_2$.  

Lemma~\ref{lem1}, applied to $J$, implies
$$
l(R/{\rm in}(\cc)) \geq 2\m(J)^2 - 2\m(J)(b_1 + b_2).
$$
If we put $\alpha=\m(J)-\m(I)$, then $\alpha\geq 0$, and
$$
l(\O_{X,Z}/\a) \geq 2(\m(I) + \aa)^2 - 2(\m(I) + \aa)\mult_Z(f).
$$
Then, by expanding the right hand side and noting that
$\m(I) \geq \frac 12 \mult_Z(f)$, we get
$$
l(\O_{X,Z}/\a) \geq 2\m(I)^2 - 2\m(I)\mult_Z(f),
$$
hence the inequality~(\ref{eq4}).
\end{proof}

We show now that Theorem~\ref{thm4} gives a positive answer to 
Question~\ref{open_question} in the codimension two case.

\begin{thm}\label{thm3}
Under the hypothesis of Theorem~\ref{thm4}, if
$e_1(I), e_2(I)$ are the Segre numbers
of $I$, then 
\begin{equation}\label{thm3-eq}
4 \m(I) e_1(I) + e_2(I) \geq 4\m(I)^2.
\end{equation}
\end{thm}

\begin{proof}
Let $W:=\Spec \O_{X,Z}$, and recall that we denote by $\M$
the closed point of $W$. 
If $I$ is a principal ideal, then the result is trivial. Henceforth,
we will assume that $I$ is not principal.
We recall now 
how the Segre numbers are defined in the two dimensional
case (see~\cite{GG} for the general case).
Consider
the fiber product diagram 
$$
\xymatrix{
V \ar[r] \ar[d] \ar[dr]^h & \Bl_I W \ar[d] \\
\Bl_{\M}W \ar[r] & W,
}
$$
where $\Bl_IW$ and $\Bl_{\M}W$ are the blow-ups along
 $I$ and $\M$, respectively.
We can write
$$
I \. \O_V = \O_V(-E_I) \quad\text{and}\quad
\M \. \O_V = \O_V(-E_{\M})
$$
for some effective Cartier divisors $E_I, E_{\M}$ on $V$.
We denote by $E_I^Z$ the union of the components of $E_I$ which are
mapped to 
$\M$, and set $E_I^{X\setminus Z} = E_I - E_I^Z$. 
The Segre numbers of $I$ are then given by
$$
e_1(I) = E_I^{X\setminus Z} \. E_{\M} 
\quad\text{and}\quad e_2(I) = - E_I \. E_I^Z.
$$
We write $I = f \. \a$,
where $\a$ is a zero dimensional ideal and $f$ defines
an effective divisor $F$ of $X$. Then we see that
$$
e_1(I) = \mult_Z(f)
$$ 
and, since $E_I = E_{\a} + h^*F$ and $E_I^Z \. h^*F = 0$,
$$
e_2(I) = - E_{\a} \. E_I^Z 
= - E_{\a}^2 + E_{\a} \. E_I^{X-Z} \geq e(\a).
$$
Therefore inequality~(\ref{thm3-eq}) follows from Theorem~\ref{thm4}.
\end{proof}

The inequality~(\ref{thm4-eq}) can be strengthened when $\dim X = 2$
by choosing carefully the local coordinates before degenerating
to monomial ideals. So, in addition
to the assumptions of Theorem~\ref{thm4},
suppose that $X$ is a surface. Then
$Z$ is a point of $X$, which we denote by $p$.
We write the ideal $I \subset \O_{X,p}$ in the form
$$
I = f \. \a,
$$
where $f \in \O_{X,p}$ and $\a$ is zero dimensional. 
If $f$ is a unit in $\O_{X,p}$, then we already know that 
inequality~(\ref{thm4-eq}) is sharp by Theorem~\ref{cor1=}. 
Henceforth, we will assume that $\mult_p(f) > 0$.
We consider the divisor $F$ on $X$ defined by $f$.
Let $E_p \cong \P^1$ be the exceptional divisor of the blow-up
of $X$ at $p$, and write the projectivized
tangent cone $\P C_pF \subset E_p$
as a divisor $c_1 P_1 + \dots + c_r P_r$, where $P_i \in E_p$ are distinct
and $c_i > 0$. Then set
\begin{equation}\label{b_i}
b_1 = c_1 + \cdots + c_{r-1} \quad\text{and}\quad b_2 = c_r.
\end{equation}
Note that $b_1 + b_2 = \mult_p(f)$ and $b_1 = 0$ if and only if $r=1$.
Let $\m(I)$ be the inverse of the log canonical threshold of $(\O_{X,p},
I)$.

\begin{thm}\label{thm5}
With the above notation, we have the following inequality:
\begin{equation}\label{thm5-eq}
4 \m(I) \mult_p(f) - 4 b_1b_2 + e(\a) \geq 4 \m(I)^2. 
\end{equation}
\end{thm}

\begin{proof}
Again, it is enough to prove that
\begin{equation}\label{eq5}
4 \m(I) \mult_p(f) - 4 b_1b_2 + 2l(\O_{X,p}/\a) \geq 4 \m(I)^2. 
\end{equation}
We proceed as in the proof of Theorem~\ref{thm4}.
We pick a regular system of parameters $x_1,x_2$ at $p$,
such that the tangent direction to $F$ at $p$, corresponding to
$P_r$, has equation $x_1=0$.
We take a multiplicative order on the coordinates
such that $x_1<x_2$.

Using the notation in the proof of Theorem~\ref{thm4},
we have $J=x_1^{b_1}x_2^{b_2}\cdot \cdot {\rm in}(\cc)$,
where $b_1$ and $b_2$ are the ones defined in (\ref{b_i}), and ${\rm in}(\cc)$
is a zero dimensional monomial ideal.
 
Again, we have $\mult_p(f) = b_1 + b_2$, 
$l(\O_{X,p}/\a) = l(R/{\rm in}(\cc))$, and $\m(I) \leq \m(J)$.
Since $J$ is a monomial ideal, Lemma~\ref{lem1} gives
$$
l(R/\a) \geq 2(\m(J) - b_1)(\m(J) - b_2).
$$
Then, arguing as in the proof of Theorem~\ref{thm4}, we obtain
(\ref{eq5}).
\end{proof}

\begin{rmk}
We remark that, if $r \geq 3$, then there is a certain freedom
in choosing $b_1, b_2$ satisfying~(\ref{b_i}) by reordering
the points $P_1,\dots,P_r$. One can check that the strongest inequality
is obtained when $c_r = \max_i \{ c_i \}$.
\end{rmk}

\begin{rmk}
In the special case when $\Supp(F)$ has exactly two smooth branches
meeting transversally at $p$, Theorem~\ref{thm5}
implies an earlier result of Corti~\cite{corti}.
\end{rmk}
 
We end with a comment about the boundary cases
in inequality~(\ref{thm5-eq}). By suitably adapting the arguments
in the proof of 
Theorem~\ref{cor1=}, we can show the following characterization:
under the assumptions of Theorem~\ref{thm5}, if
\begin{equation}\label{thm5=-eq}
4 \m(I) \mult_p(f) - 4 b_1b_2 + e(\a) = 4 \m(I)^2,
\end{equation}
then $2 \m(I) \in \N$ and, if $J$ is
as in the proof of Theorem~\ref{thm5}, then $J$ has integral closure 
$$
\ov{J} = x_1^{b_1}x_2^{b_2} \. 
\ov{(x_1^{2\m(I)-2b_1}, x_2^{2\m(I)-2b_2})}.
$$
The example below shows that, under these assumptions, 
we can not expect to get a characterization for the integral closure of $I$,
as we did for zero dimensional ideals in Theorem~\ref{cor1=}.
Consider the following ideal in $k[x_1,x_2]$:
$$
I = f \. \a = x_2^2 \. (x_1^6, x_2^2 + x_1^2x_2).
$$ 
Then $\m(I) = 3$, $\mult_p(f) = 2$, $b_1=0$, $b_2=2$ and $e(\a) = 12$, so 
equality~(\ref{thm5=-eq}) is satisfied, but
$I \not\subset x_2^2 \. \ov{(x_1^6,x_2^2)}$.

\providecommand{\bysame}{\leavevmode \hbox \o3em
{\hrulefill}\thinspace}

\end{document}